\documentclass[12pt]{amsart}
\usepackage{amscd,amssymb,amsthm,verbatim}
\newcommand{\Diff}{\operatorname{Diff}}

\newcommand{\dvol}{\operatorname{dvol}}

\newcommand{\Ham}{\operatorname{Ham}}

\newcommand{\Id}{\operatorname{Id}}

\newcommand{\Image}{\operatorname{Im}}

\newcommand{\N}{{\mathbb N}}

\newcommand{\R}{{\mathbb R}}

\numberwithin{equation}{section}
\theoremstyle{plain}
\newtheorem{definition}[equation]{Definition}

\newtheorem{lemma}[equation]{Lemma}
\newtheorem{theorem}[equation]{Theorem}
\newtheorem{proposition}[equation]{Proposition}
\newtheorem{corollary}[equation]{Corollary}

\theoremstyle{remark}

\newtheorem{remark}[equation]{Remark}

\usepackage{graphicx}
\input{amssym.def}
\input{amssym}
\input{epsf}
\usepackage{a4wide}

\def\R{\mathbb R}
\def\N{\mathbb N}

\def\Id{{\rm Id}\,}

\def\2dr#1#2{\left. \frac{d^2}{d{#1}^2} \right |_{#2}}
\def\d2#1{\frac{d^2}{d{#1}^2}}

\def\begeq{\begin{equation} \label{}}
\def\endeq{\end{equation}}
\def\begar{\begin{eqnarray}}
\def\endar{\end{eqnarray}}
\def\begar*{\begin{eqnarray*}}
\def\endar*{\end{eqnarray*}}
\def\begal{\begin{align} \label{}}
\def\endal{\end{align}}
\def\begal*{\begin{align*}}
\def\endal*{\end{align*}}

\theoremstyle{definition}

\theoremstyle{remark}

\newtheorem*{Thm*}{Theorem}
\newtheorem*{Lem*}{Lemma}
\newtheorem*{Conj*}{Conjecture}
\newtheorem*{Cor*}{Corollary}
\newtheorem*{Def*}{Definition}
\newtheorem*{Prop*}{Proposition}
\newtheorem*{Exo*}{Exercise}
\newtheorem*{Exs*}{Examples}
\newtheorem*{Ex*}{Example}
\newtheorem*{Rk*}{Remark}
\newtheorem*{Rks*}{Remarks}

\begin{document}

\title[Some geometric calculations on Wasserstein space]
{Some geometric calculations on Wasserstein space}

\author{John Lott}
\address{Department of Mathematics\\
University of Michigan\\
Ann Arbor, MI  48109-1109\\
USA} \email{lott@umich.edu}

\thanks{This research was partially 
supported by NSF grant DMS-0604829}
\date{March 31, 2007}

\begin{abstract}
We compute the Riemannian connection and curvature for the
Wasserstein space of a smooth compact Riemannian manifold.
\end{abstract}

\maketitle
\section{Introduction}

If $M$ is a smooth compact Riemannian manifold then the Wasserstein space $P_2(M)$ is the
space of Borel probability measures on $M$, equipped with the
Wasserstein metric $W_2$. We refer to \cite{Villani} for background information on 
Wasserstein spaces. The Wasserstein space originated in the study of
optimal transport.  It has had applications to PDE theory \cite{Otto},
metric geometry \cite{Lott-Villani,Sturm,Sturm2} and
functional inequalities \cite{Lott-Villani2,Otto-Villani (2000)}.

Otto showed that the heat flow on measures can be considered as a 
gradient flow on Wasserstein space \cite{Otto}.  In order to do this, he
introduced a certain formal Riemannian metric on the Wasserstein space.
This Riemannian metric has some remarkable properties.  
Using O'Neill's theorem, Otto gave
a formal argument that 
$P_2(\R^n)$ has nonnegative sectional curvature.
This was made rigorous in \cite[Theorem A.8]{Lott-Villani} and
\cite[Proposition 2.10]{Sturm} in the following sense : $M$ has nonnegative sectional
curvature if and only if the length space
$P_2(M)$ has nonnegative Alexandrov curvature. 

In this paper we study the Riemannian geometry of the Wasserstein space.
In order to write meaningful expressions, we restrict ourselves to the subspace
$P^\infty(M)$ of absolutely continuous measures with a smooth positive density function.
The space $P^\infty(M)$  is a smooth infinite-dimensional manifold in the sense,
for example, of \cite{Kriegl-Michor}.
The formal calculations that we perform can be considered as rigorous calculations
on this smooth manifold, although we do not emphasize this point.

In Section \ref{section3} we show that if $c$ is a smooth immersed
curve in $P^\infty(M)$ then its
length in $P_2(M)$, in the sense of metric geometry, equals its Riemannian
length as computed with Otto's metric.  In Section \ref{section4}
we compute the Levi-Civita connection on $P^\infty(M)$. We use it to derive the
equation for parallel transport and the geodesic equation. 

In Section \ref{section5} we compute the Riemannian curvature of $P^\infty(M)$. The
answer is relatively simple. As an application, 
if $M$ has sectional curvatures bounded below
by $r \in \R$, one can ask whether $P^\infty(M)$ necessarily has sectional curvatures
bounded below by $r$. This turns out to be the case if and only if $r = 0$.

There has been recent interest in doing Hamiltonian mechanics on the Wasserstein space of
a symplectic manifold 
\cite{Ambrosio-Gangbo,Gangbo-Nguyen-Tudorascu,Gangbo-Pacini}.  In Section \ref{section6} 
we briefly describe the Poisson geometry of $P^\infty(M)$. We show
that if $M$ is a Poisson manifold then $P^\infty(M)$ has a natural Poisson structure. We 
also show that if $M$ is symplectic then 
the symplectic leaves of the Poisson structure on $P^\infty(M)$ are the orbits of the group of
Hamiltonian diffeomorphisms, thereby making contact with \cite{Ambrosio-Gangbo,Gangbo-Pacini}.
This approach is not really new; closely related results, with applications to PDEs, 
were obtained quite a while ago by Alan Weinstein and
collaborators \cite{Marsden-Weinstein,Marsden-Ratiu-Schmid-Spencer-Weinstein,Weinstein}. 
However, it may be worth advertising this viewpoint.

I thank Wilfrid Gangbo, Tommaso Pacini and Alan Weinstein for telling me of their work.
I thank C\'edric Villani for helpful discussions and the referee for
helpful remarks.

\section{Manifolds of measures} \label{section2}

In what follows, we use the Einstein summation convention freely.

Let $M$ be a smooth connected closed Riemannian manifold of positive dimension. 
We denote the Riemannian
density by $\dvol_M$. Let $P_2(M)$ denote the space of Borel probability
measures on $M$, equipped with the Wasserstein metric $W_2$. For relevant
results about optimal transport and 
the Wasserstein metric, we refer to \cite[Sections 1 and 2]{Lott-Villani}
and references therein.

Put 
\begin{equation} \label{2.1}
P^\infty(M) \: = \: \{ \rho \: \dvol_M \: : \: \rho \in C^\infty(M), \rho > 0, \int_M \rho \:
\dvol_M \: = \: 1\}.
\end{equation}
Then $P^\infty(M)$ is a dense subset of $P_2(M)$, 
as is the complement of
$P^\infty(M)$ in $P_2(M)$.
We do not claim that $P^\infty(M)$
is necessarily a totally convex subset of $P_2(M)$, i.e. that if $\mu_0, \mu_1 \in 
P^\infty(M)$ then the minimizing geodesic in $P_2(M)$
joining them necessarily lies in $P^\infty(M)$. However,
the absolutely continuous probability measures on $M$ do form a totally convex 
subset  of $P_2(M)$ \cite{McCann}.
For the purposes of this paper,
we give $P^\infty(M)$ the smooth topology. (This differs from the subspace topology on
$P^\infty(M)$ coming from its inclusion in $P_2(M)$.) Then $P^\infty(M)$ has the structure of an
infinite-dimensional smooth manifold in the sense of
\cite{Kriegl-Michor}. The formal calculations in this paper can be rigorously 
justified as being calculations on the smooth manifold $P^\infty(M)$.  However, we 
will not belabor this point.
 
Given $\phi \in C^\infty(M)$, define $F_\phi \in C^\infty(P^\infty(M))$ by
\begin{equation} \label{2.2}
F_\phi(\rho \: \dvol_M) \:  = \: \int_M \phi \: \rho \: \dvol_M.
\end{equation} 
This gives an injection $P^\infty(M) \rightarrow (C^\infty(M))^*$,
i.e. the functions $F_\phi$ separate points in $P^\infty(M)$.
We will think of the functions $F_\phi$ as ``coordinates'' on
$P^\infty(M)$.

Given $\phi \in C^\infty(M)$, define a vector field $V_\phi$ on
$P^\infty(M)$ by saying that for $F \in C^\infty(P^\infty(M))$,
\begin{equation} \label{2.3}
(V_\phi F)(\rho \dvol_M) \: = \: 
\frac{d}{d\epsilon} \Big|_{\epsilon = 0} 
F \left( \rho \dvol_M \: - \: \epsilon \: \nabla^i ( \rho \nabla_i \phi) \dvol_M \right).
\end{equation}
The map $\phi \rightarrow V_\phi$ passes to an isomorphism
$C^\infty(M)/\R \rightarrow T_{\rho \dvol_M} P^\infty(M)$.
This parametrization of $T_{\rho \dvol_M} P^\infty(M)$ goes back to
Otto's paper
\cite{Otto}; see \cite{Ambrosio-Gigli-Savare (2004)}
for further discussion.
Otto's Riemannian metric on $P^\infty(M)$ is given \cite{Otto} by
\begin{align} \label{2.4}
\langle V_{\phi_1}, V_{\phi_2} \rangle (\rho \dvol_M) \: & = \:
\int_M \langle \nabla \phi_1, \nabla \phi_2 \rangle \: \rho \: \dvol_M \\
&  = \: - \: 
\int_M  \phi_1 \nabla^i ( \rho \nabla_i \phi_2) \: \dvol_M. \notag
\end{align}
In view of (\ref{2.3}), we write $\delta_{V_{\phi}} \rho \: = \: - \:  \nabla^i ( \rho \nabla_i \phi)$.
Then
\begin{equation} \label{2.5}
\langle V_{\phi_1}, V_{\phi_2} \rangle (\rho \dvol_M) \:  = \:
\int_M  \phi_1 \: \delta_{V_{\phi_2}} \rho \: \dvol_M
\: = \: \int_M  \phi_2 \: \delta_{V_{\phi_1}} \rho \: \dvol_M.
\end{equation}

In terms of the weighted $L^2$-spaces $L^2(M, \rho \:\dvol_M)$ and 
$\Omega^1_{L^2}(M, \rho \: \dvol_M)$, let $d$ be the usual
differential on functions and let $d_\rho^*$ be its formal adjoint.
Then (\ref{2.4}) can be written as
\begin{equation} \label{2.6}
\langle V_{\phi_1}, V_{\phi_2} \rangle (\rho \dvol_M) \: = \:
\int_M \langle d \phi_1, d \phi_2 \rangle \: \rho \: \dvol_M \:
  = \:  
\int_M  \phi_1 \: d_\rho^* d\phi_2 \: \rho \: \dvol_M. 
\end{equation}

We now relate the function $F_\phi$ and the vector field $V_\phi$.

\begin{lemma} \label{2.7}
The gradient of $F_\phi$ is $V_\phi$.
\end{lemma}
\begin{proof}
Letting $\overline{\nabla} F_\phi$ denote the gradient of $F_\phi$, 
for all $\phi^\prime \in C^\infty(M)$ we have
\begin{align} \label{2.8}
\langle \overline{\nabla} F_\phi, V_{\phi^\prime} \rangle (\rho \: \dvol_M) \: & = \:
(V_{\phi^\prime} F_\phi) (\rho \: \dvol_M) \: = \: - \: \int_M \: \phi \: 
\nabla^i ( \rho \nabla_i \phi^\prime) \dvol_M \\
& = \:
\langle V_{\phi}, V_{\phi^\prime} \rangle (\rho \: \dvol_M). \notag
\end{align}
This proves the lemma.
\end{proof}

\section{Lengths of curves} \label{section3}

In this section we relate the Riemannian metric (\ref{2.4}) to the Wasserstein metric.
One such relation was given in \cite{Otto-Villani (2000)}, where it was heuristically
shown that the geodesic distance coming from (\ref{2.4})
equals the Wasserstein metric.  To give a rigorous relation,
we recall that a curve 
$c \: : [0,1] \rightarrow P_2(M)$ has a length given by 
\begin{equation} \label{3.1}
L(c) \: = 
\: \sup_{J\in \N}\; \sup_{0 = t_0 \le t_1 \le \ldots \le t_J = 1}
\sum_{j=1}^{J} W_2\bigl(c(t_{j-1}), c(t_{j})\bigr).
\end{equation}
From the triangle inequality, the expression $\sum_{j=1}^{J} W_2\bigl(c(t_{j-1}), c(t_{j})\bigr)$
is nondecreasing under a refinement of the partition $0 = t_0 \le t_1 \le \ldots \le t_J = 1$.

If $c \: : [0,1] \rightarrow P^\infty(M)$ is a smooth curve in $P^\infty(M)$ then 
we write $c(t) \: = \: \rho(t) \: \dvol_M$ and let $\phi(t)$ satisfy
$\frac{\partial \rho}{dt} \: = \: - \: \nabla^i \left( \rho \nabla_i \phi \right)$,
where we normalize $\phi$ by requiring for example that
$\int_M \phi \: \rho \: \dvol_M \: = \: 0$. 
If $c$ is immersed then $\nabla \phi(t) \neq 0$.
The Riemannian length of $c$, as
computed using (\ref{2.4}), is
\begin{equation} \label{3.2}
\int_0^1 \langle c^\prime(t), c^\prime(t) \rangle^{\frac12} \: dt \: = \:
\int_0^1 \left( \int_M |\nabla \phi(t)|^2(m) \: \rho(t) \: \dvol_M \right)^{\frac12} \: dt.
\end{equation}
The next proposition says that this equals the length of $c$ in the metric sense.

\begin{proposition} \label{3.3}
If $c \: : [0,1] \rightarrow P^\infty(M)$ is a smooth 
immersed curve then its length $L(c)$ in
the Wasserstein space $P_2(M)$ satisfies
\begin{equation} \label{3.4}
L(c) \: = \: \int_0^1 \langle c^\prime(t), c^\prime(t) \rangle^{\frac12} \: dt.
\end{equation}
\end{proposition}
\begin{proof}
We can parametrize $c$ so that
$\int_M |\nabla \phi(t)|^2 \: \rho(t) \: \dvol_M$ is a constant $C>0$ with respect to $t$.

Let $\{S_t\}_{t \in [0,1]}$ be the one-parameter family of diffeomorphisms of $M$ given by
\begin{equation}  \label{3.5}
\frac{\partial S_t(m)}{\partial t} \: = \: (\nabla \phi(t))(S_t(m))
\end{equation}
with $S_0(m) = m$.
Then $c(t) \: = \: (S_t)_* (\rho(0) \: \dvol_M)$.

Given a partition $0 = t_0 \le t_1 \le \ldots \le t_J = 1$ of $[0,1]$,
a particular transference plan from $c(t_{j-1})$ to  $c(t_{j})$ comes from the Monge transport
$S_{t_j} \circ S_{t_{j-1}}^{-1}$. Then
\begin{align} \label{3.6}
W_2\bigl(c(t_{j-1}), c(t_{j}))^2 \: 
& \le \:
\int_M d(m, S_{t_j} (S_{t_{j-1}}^{-1}(m)))^2 \: \rho(t_{j-1}) \: \dvol_M \\
& = \:
\int_M d(S_{t_{j-1}}(m), S_{t_j}(m))^2 \: \rho(0) \: \dvol_M \notag \\
& \le \:
\int_M \left( \int_{t_{j-1}}^{t_j} |\nabla \phi(t)|(S_t(m)) \: dt \right)^2
 \: \rho(0) \: \dvol_M \notag \\
& \le \: (t_j - t_{j-1}) \: 
\int_M \int_{t_{j-1}}^{t_j} |\nabla \phi(t)|^2(S_t(m)) \: dt 
 \: \rho(0) \: \dvol_M \notag \\
& = \: (t_j - t_{j-1}) \: \int_{t_{j-1}}^{t_j}
\int_M  |\nabla \phi(t)|^2(m) 
 \: \rho(t) \: \dvol_M \: dt, \notag
\end{align}
so
\begin{align} \label{3.7}
W_2\bigl(c(t_{j-1}), c(t_{j})) \: & \le \:
(t_j - t_{j-1})^{\frac12} \: \left( \int_{t_{j-1}}^{t_j}
\int_M  |\nabla \phi(t)|^2(m) 
 \: \rho(t) \: \dvol_M \: dt \right)^{\frac12} \\
 & = \:
(t_j - t_{j-1})  \: \left(
\int_M  |\nabla \phi(t_j^\prime)|^2(m) 
 \: \rho(t_j^\prime) \: \dvol_M \right)^{\frac12} \notag
\end{align}
for some $t_j^\prime \in [t_{j-1}, t_j]$. It follows that
\begin{equation} \label{3.8}
L(c) \: \le \: \int_0^1 \langle c^\prime(t), c^\prime(t) \rangle^{\frac12} \: dt.
\end{equation}

Next, 
from \cite[Lemma A.1]{Lott-Villani},
\begin{align} \label{3.9}
& (t_j - t_{j-1}) \: \left| \int_M \phi(t_{j-1}) \: \rho(t_j) \: \dvol_M \: - \: 
\int_M \phi(t_{j-1}) \: \rho(t_{j-1}) \: \dvol_M \right|^2 \: \le \\
& W_2(c(t_{j-1}), c(t_j))^2 \: \int_{t_{j-1}}^{t_j} \int_M 
|\nabla \phi(t_{j-1})|^2 \: d\mu_t \: dt, \notag
\end{align}
where $\{\mu_t\}_{t \in [t_{j-1}, t_j]}$ is the Wasserstein geodesic
between $c(t_{j-1})$ and $c(t_j)$.
Now
\begin{align} \label{3.10}
& \int_M \phi(t_{j-1}) \: \rho(t_j) \: \dvol_M \: - \: 
\int_M \phi(t_{j-1}) \: \rho(t_{j-1}) \: \dvol_M \: = \\
& - \:
\int_M \int_{t_{j-1}}^{t_j} \phi(t_{j-1}) \:
\nabla^i \left( \rho(t) \nabla_i \phi(t) \right) \: dt \: \dvol_M \: = 
\notag \\
&\int_{t_{j-1}}^{t_j}  \int_M \langle \nabla \phi(t_{j-1}),
\nabla \phi(t) \rangle \: \rho(t) \: \dvol_M \: dt, \notag
\end{align}
so (\ref{3.9}) becomes
\begin{align} \label{3.11}
& (t_j - t_{j-1}) \:
\left( \int_{t_{j-1}}^{t_j}  \int_M \langle \nabla \phi(t_{j-1}),
\nabla \phi(t) \rangle \: \rho(t) \: \dvol_M \: dt \right)^2 \: \le \\
& W_2(c(t_{j-1}), c(t_j))^2 \: \int_{t_{j-1}}^{t_j} \int_M 
|\nabla \phi(t_{j-1})|^2 \: d\mu_t \: dt. \notag
\end{align}
Thus
\begin{equation} \label{3.12}
L(c) \: \ge \:
\sum_{j=1}^J 
\frac{
\frac{
\int_{t_{j-1}}^{t_j}  \int_M \langle \nabla \phi(t_{j-1}),
\nabla \phi(t) \rangle \: \rho(t) \: \dvol_M \: dt }{t_j - t_{j-1}}
}{\sqrt{ \frac{1}{t_j - t_{j-1}} \int_{t_{j-1}}^{t_j}
\int_M 
|\nabla \phi(t_{j-1})|^2 \: 
d\mu_t \: dt}}
 \: (t_j - t_{j-1}). 
\end{equation}

As the partition of $[0,1]$ becomes finer, the term $\frac{
\int_{t_{j-1}}^{t_j}  \int_M \langle \nabla \phi(t_{j-1}),
\nabla \phi(t) \rangle \: \rho(t) \: \dvol_M \: dt }{t_j - t_{j-1}}$ 
uniformly approaches the constant $C$. 

The Wasserstein geodesic $\{\mu_t\}_{t \in [t_{j-1}, t_j]}$
has the form $\mu_t \: = \: (F_t)_* \mu_{t_{j-1}}$ for measurable maps 
$F_t \: : \: M \rightarrow M$ with $F_{t_{j-1}} = \Id$
\cite{McCann}. Then
\begin{align}
& \left| \frac{1}{t_j - t_{j-1}} \int_{t_{j-1}}^{t_j}
\int_M 
|\nabla \phi(t_{j-1})|^2 \: 
d\mu_t \: dt \: - \: C \right| \: = \\
& \left| \frac{1}{t_j - t_{j-1}} \int_{t_{j-1}}^{t_j}
\left(
\int_M 
|\nabla \phi(t_{j-1})|^2 \: 
d\mu_t \: - \: \int_M 
|\nabla \phi(t_{j-1})|^2 \: 
d\mu_{t_{j-1}} \right) \: dt \right| \: = \notag \\
& \left| \frac{1}{t_j - t_{j-1}} \int_{t_{j-1}}^{t_j}
\int_M \left(
|\nabla \phi(t_{j-1})|^2 \circ F_t \: - \: 
|\nabla \phi(t_{j-1})|^2 \: \right) \:
d\mu_{t_{j-1}}  \: dt \right| \: \le \notag \\
& \frac{1}{t_j - t_{j-1}} \parallel \nabla |\nabla \phi(t_{j-1})|^2 
\parallel_\infty \: \int_{t_{j-1}}^{t_j} \: \int_M d(m, F_t(m)) \: 
d\mu_{t_{j-1}}(m) \: dt \: \le \notag \\
& \frac{1}{t_j - t_{j-1}} \parallel \nabla |\nabla \phi(t_{j-1})|^2 
\parallel_\infty \: \int_{t_{j-1}}^{t_j} \: \sqrt{ \int_M d(m, F_t(m))^2 \: 
d\mu_{t_{j-1}}(m)}\: dt \: = \notag \\
& \frac{1}{t_j - t_{j-1}} \parallel \nabla |\nabla \phi(t_{j-1})|^2 
\parallel_\infty \: \int_{t_{j-1}}^{t_j} \: W_2(\mu_{t_{j-1}}, \mu_t) \: dt \: \le \notag \\
& \parallel \nabla |\nabla \phi(t_{j-1})|^2 
\parallel_\infty \: W_2(c(t_{j-1}), c(t_j)). \notag
\end{align}
Now continuity of a $1$-parameter 
family of smooth measures in the smooth topology implies 
continuity in the weak-$*$ topology, which is metricized by
$W_2$ (as $M$ is compact). It follows that
as the partition of $[0,1]$ becomes finer, the term  
$\frac{1}{t_j - t_{j-1}} \int_{t_{j-1}}^{t_j}
\int_M 
|\nabla \phi(t_{j-1})|^2 \: 
d\mu_t \: dt$
uniformly approaches the constant $C$. Thus from (\ref{3.12}),
\begin{equation} \label{3.13}
L(c) \: \ge \: \sqrt{C} \: = \:  \int_0^1 \langle c^\prime(t), c^\prime(t) \rangle^{\frac12} \: dt.
\end{equation}
This proves the proposition.
\end{proof}

\begin{remark} \label{3.14}
Let $X$ be a finite-dimensional Alexandrov space and let $R$ be its set of
nonsingular points.  There is a continuous Riemannian metric $g$ on
$R$ so that lengths of curves in $R$ can be computed using $g$
\cite{Otsu-Shioya}. (Note that in general, $R$ and $X-R$ are
dense in $X$.)  
This is somewhat similar to the situation for
$P^\infty(M) \subset P_2(M)$. 

In fact, there is an open dense
subset $O \subset X$ with a Lipschitz manifold structure and a Riemannian
metric of bounded variation that extends $g$ \cite{Perelman}.
We do not know if there is a Riemannian manifold structure, in some
appropriate sense, on an open dense subset of $P_2(M)$. 
Other approaches to geometrizing $P_2(M)$, with a view toward
gradient flow, are in
\cite{Ambrosio-Gigli-Savare (2004),Carrillo-McCann-Villani (2004)};
see also \cite{Ohta}.
\end{remark}

\section{Levi-Civita connection, parallel transport and geodesics} \label{section4}

In this section we compute the Levi-Civita connection of $P^\infty(M)$.
We derive the formula for parallel transport in $P^\infty(M)$ and
the geodesic equation for $P^\infty(M)$.

We first compute commutators of our canonical vector fields $\{ V_\phi \}_{\phi \in C^\infty(M)}$.

\begin{lemma} \label{4.1}
Given $\phi_1, \phi_2 \in C^\infty(M)$, the commutator
$[V_{\phi_1}, V_{\phi_2}]$ is given by
\begin{align} \label{4.2}
& \left( [V_{\phi_1}, V_{\phi_2}]F \right) (\rho \dvol_M) \: = \\
& \frac{d}{d\epsilon} \Big|_{\epsilon = 0} F \left(
\rho \dvol_M \: - \: \epsilon \nabla_i \left[ \rho \left( (
\nabla^i \nabla^j \phi_2 ) \nabla_j \phi_1 \: - \: 
(\nabla^i \nabla^j \phi_1) \nabla_j \phi_2 \right) \right] \dvol_M
\right) \notag
\end{align}
for  $F \in C^\infty(P^\infty(M))$.
\end{lemma}
\begin{proof}
We have
\begin{align} \label{4.3}
& \left( [V_{\phi_1}, V_{\phi_2}]F \right) (\rho \dvol_M) \: = \:
\left( V_{\phi_1} (V_{\phi_2}F) \right) (\rho \dvol_M) \: - \:
\left( V_{\phi_2} (V_{\phi_1}F) \right) (\rho \dvol_M) \: = \: \\
& \frac{d}{d\epsilon_1}\Big|_{\epsilon_1 = 0} 
(V_{\phi_2}F) \left( \rho \dvol_M \: - \: \epsilon_1 \: \nabla^i ( \rho \nabla_i \phi_1) \dvol_M \right)
\: - \notag \\
& \frac{d}{d\epsilon_2} \Big|_{\epsilon_2 = 0} 
(V_{\phi_1}F) \left( \rho \dvol_M \: - \: \epsilon_2 \: \nabla^i ( \rho \nabla_i \phi_2) \dvol_M \right)
\: = \: \notag \\
& \frac{d}{d\epsilon_1} \Big|_{\epsilon_1 = 0} \:
\frac{d}{d\epsilon_2} \Big|_{\epsilon_2 = 0} \:
F \left( (\rho \: - \: \epsilon_1 \: \nabla^i ( \rho \nabla_i \phi_1)) \dvol_M \: 
- \: \epsilon_2 \: \nabla^j ( (\rho - \epsilon_1 \nabla^i(\rho \nabla_i \phi_1)) \nabla_j \phi_2) \dvol_M 
\right) \: - \: \notag \\
& \frac{d}{d\epsilon_2} \Big|_{\epsilon_2 = 0} \:
\frac{d}{d\epsilon_1} \Big|_{\epsilon_1 = 0} \:
F \left( (\rho \: - \: \epsilon_2 \: \nabla^i ( \rho \nabla_i \phi_2)) \dvol_M \: 
- \: \epsilon_1 \: \nabla^j ( (\rho - \epsilon_2 \nabla^i(\rho \nabla_i \phi_2)) \nabla_j \phi_1) \dvol_M 
\right) \: = \: \notag \\
& \frac{d}{d\epsilon} \Big|_{\epsilon = 0} \:
F \left( \rho \: \dvol_M \: + \: \epsilon \: \nabla^j (\nabla^i(\rho \nabla_i \phi_1) \nabla_j \phi_2) \dvol_M 
\: - \: \epsilon \: \nabla^j (\nabla^i(\rho \nabla_i \phi_2) \nabla_j \phi_1) \dvol_M. \notag
\right)
\end{align}
One can check that
\begin{equation} \label{4.4}
\nabla^j (\nabla^i(\rho \nabla_i \phi_1) \nabla_j \phi_2)
\: - \: \nabla^j (\nabla^i(\rho \nabla_i \phi_2) \nabla_j \phi_1) \: = \:
- \: \nabla_i \left[ \rho \left( (
\nabla^i \nabla^j \phi_2 ) \nabla_j \phi_1 \: - \: 
(\nabla^i \nabla^j \phi_1) \nabla_j \phi_2 \right) \right],
\end{equation}
from which the lemma follows.
\end{proof}

We now compute the Levi-Civita connection.

\begin{proposition} \label{4.5}
The Levi-Civita connection $\overline{\nabla}$ of $P^\infty(M)$ is given by
\begin{equation} \label{4.6}
((\overline{\nabla}_{V_{\phi_1}} V_{\phi_2}) F)(\rho \: \dvol_M) \: = \:
\frac{d}{d\epsilon}\Big|_{\epsilon = 0}
F \left( \rho \: \dvol_M \: - \: \epsilon \: \nabla_i \left( \rho \: \nabla_j \phi_1 \:
\nabla^i \nabla^j \phi_2 \right) \: \dvol_M \right)
\end{equation}
for $F \in C^\infty(P^\infty(M))$.
\end{proposition}
\begin{proof}
Define a vector field $D_{V_{\phi_1}} V_{\phi_2}$ by
\begin{equation} \label{4.7}
((D_{V_{\phi_1}} V_{\phi_2}) F)(\rho \: \dvol_M) \: = \:
\frac{d}{d\epsilon}\Big|_{\epsilon = 0}
F \left( \rho \: \dvol_M \: - \: \epsilon \: \nabla_i \left( \rho \: \nabla_j \phi_1 \:
\nabla^i \nabla^j \phi_2 \right) \: \dvol_M \right)
\end{equation}
for $F \in C^\infty(P^\infty(M))$.
We also write
\begin{equation} \label{4.8}
\delta_{D_{V_{\phi_1}} V_{\phi_2}} \rho \: = \: - \:  \nabla_i \left( \rho \: \nabla_j \phi_1 \:
\nabla^i \nabla^j \phi_2 \right).
\end{equation}
It is clear from Lemma \ref{4.1} that 
\begin{equation} \label{4.9}
D_{V_{\phi_1}} V_{\phi_2} \: - \: D_{V_{\phi_2}} V_{\phi_1}
\: = \: [V_{\phi_1}, V_{\phi_2}].
\end{equation}
Next, 
\begin{align} \label{4.10}
\left( V_{\phi_1} \langle V_{\phi_2}, V_{\phi_3} \rangle \right)(\rho \: \dvol_M) \: = \:
& - \int_M \nabla^i \phi_2 \: \nabla_i \phi_3 \: \nabla^j(\rho \nabla_j \phi_1) \: \dvol_M \\
= \: &
\int_M \nabla_j \phi_1 \: \nabla^i \nabla^j \phi_2 \: \nabla_i \phi_3 \: \rho \: \dvol_M \: + \:
\int_M \nabla_j \phi_1 \: \nabla^i \nabla^j \phi_3 \: \nabla_i \phi_2 \: \rho \: \dvol_M \notag \\
= \: & - \: 
\int_M \phi_3 \: \nabla_i (\rho \: \nabla_j \phi_1 \: \nabla^i \nabla^j \phi_2) \: \dvol_M \: - 
\notag \\
&
\int_M \phi_2 \: \nabla_i (\rho \: \nabla_j \phi_1 \: \nabla^i \nabla^j \phi_3) \: \dvol_M \notag \\
= \: & \: \int_M \phi_3 \: \delta_{D_{V_{\phi_1}} V_{\phi_2}} \rho \: \dvol_M \: + \:
\int_M \phi_2 \: \delta_{D_{V_{\phi_1}} V_{\phi_3}} \rho \: \dvol_M \notag \\
= \: & \: \langle D_{V_{\phi_1}} V_{\phi_2}, V_{\phi_3}  \rangle (\rho \: \dvol_M)\: + \:
\langle V_{\phi_2}, D_{V_{\phi_1}} V_{\phi_3} \rangle (\rho \: \dvol_M). \notag
\end{align}
Thus 
\begin{equation} \label{4.11}
V_{\phi_1} \langle V_{\phi_2}, V_{\phi_3} \rangle \: = \: 
\langle D_{V_{\phi_1}} V_{\phi_2}, V_{\phi_3}  \rangle \: + \:
\langle V_{\phi_2}, D_{V_{\phi_1}} V_{\phi_3} \rangle.
\end{equation}
As
\begin{align} \label{4.12}
2 \langle \overline{\nabla}_{V_{\phi_1}} V_{\phi_2}, V_{\phi_3} \rangle \: = \: &
V_{\phi_1} \langle V_{\phi_2}, V_{\phi_3} \rangle \: + \:
V_{\phi_2} \langle V_{\phi_3}, V_{\phi_1} \rangle \: - \:
V_{\phi_3} \langle V_{\phi_1}, V_{\phi_2} \rangle \: + \\
& \langle V_{\phi_3}, [V_{\phi_1}, V_{\phi_2}] \rangle \: - \: 
\langle V_{\phi_2}, [V_{\phi_1}, V_{\phi_3}] \rangle \: - \: 
\langle V_{\phi_1}, [V_{\phi_2}, V_{\phi_3}] \rangle, \notag
\end{align}
substituting (\ref{4.9}) and (\ref{4.11}) into the right-hand side of (\ref{4.12}) shows that
\begin{equation} \label{4.13}
\langle \overline{\nabla}_{V_{\phi_1}} V_{\phi_2}, V_{\phi_3} \rangle \: = \: 
\langle D_{V_{\phi_1}} V_{\phi_2}, V_{\phi_3} \rangle
\end{equation}
for all $\phi_3 \in C^\infty(M)$. The proposition follows.
\end{proof}

\begin{lemma} \label{4.14}
The connection coefficients at $\rho \: \dvol_M$ are given by
\begin{equation} \label{4.15}
\langle \overline{\nabla}_{V_{\phi_1}} V_{\phi_2}, V_{\phi_3} \rangle \: = \:
\int_M \nabla_i \phi_1 \: \nabla_j \phi_3 \: \nabla^i \nabla^j \phi_2 \: \rho \: \dvol_M.
\end{equation}
\end{lemma}
\begin{proof}
This follows from (\ref{2.5}) and (\ref{4.6}).
\end{proof}

Let $G_\rho$ be the Green's operator for $d_\rho^* d$ on $L^2(M, \rho \: \dvol_M)$. (More explicitly, if
$\int_M f \: \rho \: \dvol_M \: = \: 0$ and $\phi \: = \:
G_\rho f$ then $\phi$ satisfies
$- \: \frac{1}{\rho} \:
\nabla^i (\rho \nabla_i \phi) \: = \: f$ and
$\int_M \phi \: \: \rho \dvol_M \: = \: 0$, while $G_\rho 1 \: = \: 0$.)
Let $\Pi_\rho$ denote orthogonal projection onto $\overline{\Image(d)} \: \subset \:
\Omega^1_{L^2}(M, \rho \: \dvol_M)$.

\begin{lemma} \label{4.16}
At $\rho \: \dvol_M$, we have
$\overline{\nabla}_{V_{\phi_1}} V_{\phi_2} \: = \: V_{\phi}$, where
$\phi \: = \: G_\rho d_\rho^* (\nabla_i \nabla_j \phi_2 \: \nabla^j \phi_1 \: dx^i)$.
\end{lemma}
\begin{proof}
Given $\phi_3 \in C^\infty(M)$, we have
\begin{align} \label{4.17}
\langle V_{\phi_3}, V_{\phi} \rangle (\rho \: \dvol_M) \: & = \:
\int_M \langle d\phi_3, d G_\rho d_\rho^* (\nabla_i \nabla_j \phi_2 \: \nabla^j \phi_1 \: dx^i)
\rangle \: \rho \: \dvol_M \\
& = \: \int_M \langle d\phi_3, 
\Pi_\rho (\nabla_i \nabla_j \phi_2 \: \nabla^j \phi_1 \: dx^i)
\rangle \: \rho \: \dvol_M \notag \\
& = \:  \int_M \langle d\phi_3, 
\nabla_i \nabla_j \phi_2 \: \nabla^j \phi_1 \: dx^i
\rangle \: \rho \: \dvol_M \: = \: \langle V_{\phi_3}, \overline{\nabla}_{V_{\phi_1}} V_{\phi_2}
\rangle (\rho \: \dvol_M).
\notag
\end{align}
The lemma follows.
\end{proof}

To derive the equation for parallel transport, let
$c \: : \: (a,b) \rightarrow P^\infty(M)$ be a smooth curve. As before, we write
$c(t) \: = \: \rho(t) \: \dvol_M$ and define $\phi(t) \in C^\infty(M)$, up to a constant, 
by $\frac{dc}{dt} \: = \: V_{\phi(t)}$. Let $V_{\eta(t)}$ be a vector field along $c$,
with $\eta(t) \in C^\infty(M)$.
If $\{\phi_\alpha\}_{\alpha=1}^\infty$ is a basis for $C^\infty(M)/\R$ then
$\{V_{\phi_\alpha} \}_{\alpha=1}^\infty$ is a global basis for $TP^\infty(M)$ and we
can write $\eta(t) \: = \: \sum_\alpha \eta_\alpha(t) \: V_{\phi_\alpha} \Big|_{c(t)}$.
The condition for $V_{\eta}$ to be parallel along $c$ is
\begin{equation} \label{4.18}
\sum_\alpha \frac{d\eta_\alpha}{dt} \: V_{\phi_\alpha} \Big|_{c(t)} \: + \: 
\sum_\alpha \eta_\alpha(t) \: \overline{\nabla}_{V_{\phi(t)}} V_{\eta_\alpha} \Big|_{c(t)} \: = \: 0,
\end{equation}
or
\begin{equation} \label{4.19}
 V_{\frac{\partial \eta}{\partial t}} \: + \: 
\overline{\nabla}_{V_{\phi(t)}} V_{\eta(t)} \: = \: 0.
\end{equation}
\begin{proposition} \label{4.20}
The equation for $V_{\eta}$ to be parallel along $c$ is
\begin{equation} \label{4.21}
\nabla_i \left( \rho \left( \nabla^i \frac{\partial \eta}{\partial t} \: + \:
\nabla_j \phi \: \nabla^i \nabla^j \eta \right) \right) \: = \: 0.
\end{equation}
\end{proposition}
\begin{proof}
This follows from (\ref{2.3}), (\ref{4.6}) and (\ref{4.19}).
\end{proof}

As a check on equation (\ref{4.21}), we show that parallel transport along $c$
preserves the inner product.

\begin{lemma} \label{4.22}
If $V_{\eta_1}$ and $V_{\eta_2}$ are parallel vector fields along $c$ then
$\int_M \langle \nabla \eta_1, \nabla \eta_2 \rangle \: \rho \: \dvol_M$ is
constant in $t$.
\end{lemma}
\begin{proof}
We have
\begin{align} \label{4.23}
\frac{d}{dt} \int_M \langle \nabla \eta_1, \nabla \eta_2 \rangle \: \rho \: \dvol_M
\: = \: &
\int_M \nabla^i \frac{\partial \eta_1}{\partial t} \: \nabla_i \eta_2 \: 
\rho \: \dvol_M \: + \:
\int_M \nabla_i \eta_1 \: \nabla^i \frac{\partial \eta_2}{dt} \: \rho \: \dvol_M  - \: \\
& \int_M \nabla_i \eta_1 \nabla^i \eta_2 \: \nabla^j (\rho \nabla_j \phi)\: \dvol_M \notag \\
\: = \: & \int_M \nabla^i \frac{\partial \eta_1}{\partial t} \: \nabla_i \eta_2 \: 
\rho \: \dvol_M \: + \:
\int_M \nabla_i \eta_1 \: \nabla^i \frac{\partial \eta_2}{dt} \: \rho \: \dvol_M  + \: \notag \\
& \int_M \left( \nabla^i \nabla^j \eta_1 \: \nabla_i \eta_2 \: + \:
\nabla_i \eta_1 \: \nabla^i \nabla^j \eta_2  \right) \: \nabla_j \phi\: \rho \: \dvol_M \notag \\
\: = \: & - \: \int_M \eta_2 \: 
\nabla_i \left( \rho \left( \nabla^i \frac{\partial \eta_1}{\partial t} \: + \:
\nabla_j \phi \: \nabla^i \nabla^j \eta_1 \right) \right) \: \dvol_M \: - \notag \\
& \int_M \eta_1 \: 
\nabla_i \left( \rho \left( \nabla^i \frac{\partial \eta_2}{\partial t} \: + \:
\nabla_j \phi \: \nabla^i \nabla^j \eta_2 \right) \right) \: \dvol_M \notag \\
\: = \: & 0. \notag
\end{align}
This proves the lemma.
\end{proof}

Finally, we derive the geodesic equation.

\begin{proposition} \label{4.24}
The geodesic equation for $c$ is 
\begin{equation} \label{4.25}
\frac{\partial \phi}{\partial t} \: + \: \frac12 \: |\nabla \phi|^2 \: = \: 0,
\end{equation}
modulo the addition of a spatially-constant function to $\phi$.
\end{proposition}
\begin{proof}
Taking $\eta \: = \: \phi$ in (\ref{4.21}) gives
\begin{equation} \label{4.26}
\nabla_i \left( \rho \: \nabla^i \left( \frac{\partial \phi}{\partial t} \: + \: \frac12 \:
|\nabla \phi|^2 \right) \right) \: = \: 0.
\end{equation}
Thus $\frac{\partial \phi}{\partial t} \: + \: \frac12 \:
|\nabla \phi|^2$ is spatially constant.  Redefining $\phi$ by adding to it a
function of $t$ alone, we can assume that $(\ref{4.25})$ holds.
\end{proof}

\begin{remark} \label{4.27}
Equation (\ref{4.25}) has been known for a while, 
at least in the case of $\R^n$, to be the formal equation for
Wasserstein geodesics. For general
Riemannian manifolds $M$, 
it was formally derived as the Wasserstein geodesic equation in
\cite{Otto-Villani (2000)} by minimizing lengths of curves. For $t > 0$, it has the Hopf-Lax solution
\begin{equation} \label{4.28}
\phi(t,m) \: = \: \inf_{m^\prime \in M} \left( \phi(0, m^\prime) \: + \: 
\frac{d(m,m^\prime)^2}{2t} \right).
\end{equation}

Given $\mu_0, \mu_1 \in P^\infty(M)$, it is known that there is a unique minimizing Wasserstein
geodesic $\{ \mu_t \}_{t \in [0,1]}$ joining them.  It is of the form
$\mu_t \: = \: (F_t)_* \mu_0$, where $F_t  \in \Diff(M)$ is given by
$F_t(m) \:  = \: \exp_m (-  t \nabla_m \phi_0)$ for an appropriate Lipschitz
function $\phi_0$ \cite{McCann}. If $\phi_0$ happens to be smooth then defining $\rho(t)$ by
$\mu_t \: = \: \rho(t) \: \dvol_M$ and defining $\phi(t) \in C^\infty(M)/\R$ as above,
it is known that $\phi$ satisfies (\ref{4.25}), with $\phi(0) \: = \: \phi_0$
\cite[Section 5.4.7]{Villani}. In this way,
(\ref{4.25}) rigorously describes certain geodesics in the Wasserstein space $P_2(M)$.
\end{remark}

\section{Curvature} \label{section5}

In this section we compute the Riemannian curvature tensor of $P^\infty(M)$.

Given $\phi, \phi^\prime \in C^\infty(M)$, define
$T_{\phi \phi^\prime} \in \Omega^1_{L^2}(M)$ by
\begin{equation} \label{5.1}
T_{\phi \phi^\prime} \: = \: (I - \Pi_\rho) \: \left( \nabla^i \phi \: \nabla_i \nabla_j \phi^\prime \:
dx^j \right).
\end{equation}
(The left-hand side depends on $\rho$, but we suppress this for simplicity of notation.)

\begin{lemma} \label{5.2}
$T_{\phi \phi^\prime} \: + \: T_{\phi^\prime \phi} \: = \: 0$.
\end{lemma}
\begin{proof}
As 
\begin{equation} \label{5.3}
\nabla^i \phi \: \nabla_i \nabla_j \phi^\prime \:
dx^j \: + \: \nabla^i \phi^\prime \: \nabla_i \nabla_j \phi \:
dx^j \: = \: d \langle \nabla \phi, \nabla \phi^\prime \rangle,
\end{equation}
and $I - \Pi_\rho$ projects away from $\Image(d)$, 
the lemma follows.
\end{proof}

\begin{theorem} \label{5.4}
Given $\phi_1, \phi_2, \phi_3, \phi_4 \in C^\infty(M)$, the Riemannian curvature
operator $\overline{R}$ of $P^\infty(M)$ is given by
\begin{align} \label{5.5}
\langle \overline{R}(V_{\phi_1}, V_{\phi_2}) V_{\phi_3}, V_{\phi_4} \rangle \: = \: &
\int_M \langle R(\nabla \phi_1, \nabla \phi_2) \nabla \phi_3, \nabla \phi_4 \rangle \: \rho \:
\dvol_M \: - \: 2 \langle T_{\phi_1 \phi_2}, T_{\phi_3 \phi_4} \rangle \: + \\
& \langle T_{\phi_2 \phi_3}, T_{\phi_1 \phi_4} \rangle \: - \:
\langle T_{\phi_1 \phi_3}, T_{\phi_2 \phi_4} \rangle, \notag
\end{align}
where both sides are evaluated at $\rho \: \dvol_M \in P^\infty(M)$.
\end{theorem}
\begin{proof}
We use the formula
\begin{align} \label{5.6}
\langle \overline{R}(V_{\phi_1}, V_{\phi_2}) V_{\phi_3}, V_{\phi_4} \rangle \: = \: &
V_{\phi_1} \langle \overline{\nabla}_{V_{\phi_2}} V_{\phi_3}, V_{\phi_4} \rangle \: - \:
\langle \overline{\nabla}_{V_{\phi_2}} V_{\phi_3}, \overline{\nabla}_{V_{\phi_1}} V_{\phi_4} \rangle \:
- \\
&  V_{\phi_2} \langle \overline{\nabla}_{V_{\phi_1}} V_{\phi_3}, V_{\phi_4} \rangle \: + \:
\langle \overline{\nabla}_{V_{\phi_1}} V_{\phi_3}, \overline{\nabla}_{V_{\phi_2}} V_{\phi_4} \rangle
\: - \notag \\
& \langle \overline{\nabla}_{[V_{\phi_1}, V_{\phi_2}]} V_{\phi_3}, V_{\phi_4} \rangle. \notag
\end{align}

First, from (\ref{2.3}) and (\ref{4.14}),
\begin{align} \label{5.7} 
V_{\phi_1} \langle \overline{\nabla}_{V_{\phi_2}} V_{\phi_3}, V_{\phi_4} \rangle \: = \: &
- \: \int_M \nabla_i \phi_2 \: \nabla_j \phi_4 \: \nabla^i \nabla^j \phi_3 \: \nabla^k(
\rho \nabla_k \phi_1) \: \dvol_M \\
= \: &
\int_M \nabla^k \nabla_i \phi_2 \: \nabla_j \phi_4 \: \nabla^i \nabla^j \phi_3 \: 
\nabla_k \phi_1 \: \rho \: \dvol_M \: + \notag \\
& \int_M \nabla_i \phi_2 \: \nabla^k \nabla_j \phi_4 \: \nabla^i \nabla^j \phi_3 \: 
\nabla_k \phi_1 \: \rho \: \dvol_M \: + \notag \\
& \int_M \nabla_i \phi_2 \: \nabla_j \phi_4 \: \nabla^k \nabla^i \nabla^j \phi_3 \: 
\nabla_k \phi_1 \: \rho \:  \dvol_M. \notag
\end{align}
Similarly,
\begin{align} \label{5.8}
V_{\phi_2} \langle \overline{\nabla}_{V_{\phi_1}} V_{\phi_3}, V_{\phi_4} \rangle \: = \:
&
\int_M \nabla^k \nabla_i \phi_1 \: \nabla_j \phi_4 \: \nabla^i \nabla^j \phi_3 \: 
\nabla_k \phi_2 \: \rho \: \dvol_M \: +  \\
& \int_M \nabla_i \phi_1 \: \nabla^k \nabla_j \phi_4 \: \nabla^i \nabla^j \phi_3 \: 
\nabla_k \phi_2 \: \rho \: \dvol_M \: + \notag \\
& \int_M \nabla_i \phi_1 \: \nabla_j \phi_4 \: \nabla^k \nabla^i \nabla^j \phi_3 \: 
\nabla_k \phi_2 \: \rho \: \dvol_M. \notag
\end{align}

Next, using (\ref{2.4}), Lemma \ref{4.16} and (\ref{5.1}),
\begin{align} \label{5.9}
\langle \overline{\nabla}_{V_{\phi_2}} V_{\phi_3}, \overline{\nabla}_{V_{\phi_1}} V_{\phi_4} \rangle
\: = \: &
\langle d G_\rho d_\rho^* (\nabla_i \nabla_j \phi_3 \: \nabla^j \phi_2 \: dx^i),
d G_\rho d_\rho^* (\nabla_k \nabla_l \phi_4 \: \nabla^l \phi_1 \: dx^k) \rangle_{L^2} \\
= \: &
\langle \Pi_\rho (\nabla_i \nabla_j \phi_3 \: \nabla^j \phi_2 \: dx^i),
\Pi_\rho  (\nabla_k \nabla_l \phi_4 \: \nabla^l \phi_1 \: dx^k) \rangle_{L^2} \notag \\
\: = \: &
\langle \nabla_i \nabla_j \phi_3 \: \nabla^j \phi_2 \: dx^i,
\nabla_k \nabla_l \phi_4 \: \nabla^l \phi_1 \: dx^k \rangle_{L^2} \: - \: 
\langle T_{\phi_2 \phi_3}, T_{\phi_1 \phi_4} \rangle \notag \\
\: = \: &
\int_M  \nabla_i \nabla_j \phi_3 \: \nabla^j \phi_2 \: 
\nabla^i \nabla_l \phi_4 \: \nabla^l \phi_1 \: \rho \: \dvol_M \: - \: 
\langle T_{\phi_2 \phi_3}, T_{\phi_1 \phi_4} \rangle. \notag
\end{align}
Similarly,
\begin{equation} \label{5.10}
\langle \overline{\nabla}_{V_{\phi_1}} V_{\phi_3}, \overline{\nabla}_{V_{\phi_2}} V_{\phi_4} \rangle
\: = \: \int_M  \nabla_i \nabla_j \phi_3 \: \nabla^j \phi_1 \: 
\nabla^i \nabla_l \phi_4 \: \nabla^l \phi_2 \: \rho \: \dvol_M \: - \: 
\langle T_{\phi_1 \phi_3}, T_{\phi_2 \phi_4} \rangle.
\end{equation}

Finally, we compute 
$\langle \overline{\nabla}_{[V_{\phi_1}, V_{\phi_2}]} V_{\phi_3}, V_{\phi_4} \rangle$.
From (\ref{4.2}), we can write $[V_{\phi_1}, V_{\phi_2}] \: = \: V_{\phi}$, where
\begin{equation} \label{5.11}
\phi \: = \: G_\rho \: d_\rho^* \: \left(  \nabla_i \nabla_j \phi_2 \: \nabla^j \phi_1 
\: dx^i \: - \: \nabla_i \nabla_j \phi_1 \: \nabla^j \phi_2
\: dx^i \right).
\end{equation}
Then from (\ref{4.15}),
\begin{align} \label{5.12}
\langle \overline{\nabla}_{[V_{\phi_1}, V_{\phi_2}]} V_{\phi_3}, V_{\phi_4} \rangle \: = \:
& \int_M \nabla_i \phi \: \nabla_j \phi_4 \: \nabla^i \nabla^j \phi_3 \: \rho \: \dvol_M \: = \:
\langle d\phi, \nabla^j \phi_4 \: \nabla_i \nabla_j \phi_3 \: dx^i \rangle_{L^2} \\
= \: &
\langle d G_\rho \: d_\rho^* \: \left(  \nabla_i \nabla_j \phi_2 \: \nabla^j \phi_1 
\: dx^i \: - \: \nabla_i \nabla_j \phi_1 \: \nabla^j \phi_2
\: dx^i \right), \nabla^j \phi_4 \: \nabla_i \nabla_j \phi_3 \: dx^i \rangle_{L^2} \notag \\
= \: &
\langle \Pi_\rho \: \left(  \nabla_i \nabla_j \phi_2 \: \nabla^j \phi_1 
\: dx^i \: - \: \nabla_i \nabla_j \phi_1 \: \nabla^j \phi_2
\: dx^i \right), \Pi_\rho \left( \nabla^j \phi_4 \: \nabla_i \nabla_j \phi_3 \: dx^i \right) 
\rangle_{L^2} \notag \\
= \: & \int_M \left( \nabla_i \nabla_j \phi_2 \: \nabla^j \phi_1 
 \: - \: \nabla_i \nabla_j \phi_1 \: \nabla^j \phi_2 \right) \:  \nabla_k \phi_4 \: \nabla^i \nabla^k \phi_3
 \: \rho \: \dvol_M \: - \notag \\
& \: \langle T_{\phi_1 \phi_2}, T_{\phi_4 \phi_3} \rangle \: + \: 
\langle T_{\phi_2 \phi_1}, T_{\phi_4 \phi_3} \rangle \notag \\
= \: & \int_M \left( \nabla_i \nabla_j \phi_2 \: \nabla^j \phi_1 
 \: - \: \nabla_i \nabla_j \phi_1 \: \nabla^j \phi_2 \right) \:  \nabla_k \phi_4 \: \nabla^i \nabla^k \phi_3
 \: \rho \: \dvol_M \: +\notag \\
& \: 2 \: \langle T_{\phi_1 \phi_2}, T_{\phi_3 \phi_4} \rangle. \notag
\end{align}
The theorem follows from combining equations (\ref{5.6})-(\ref{5.12}).
\end{proof}

\begin{corollary} \label{5.13}
Suppose that $\phi_1, \phi_2 \in C^\infty(M)$ satisfy
$\int_M |\nabla \phi_1|^2 \: \rho \: \dvol_M \: = \: \int_M |\nabla \phi_2|^2 \: \rho \: \dvol_M \: = \: 1$
and $\int_M \langle \nabla \phi_1, \nabla \phi_2 \rangle \: \rho \: \dvol_M \: = \: 0$.
Then the sectional curvature at $\rho \: \dvol_M \in P^\infty(M)$ of the $2$-plane
spanned by $V_{\phi_1}$ and $V_{\phi_2}$ is
\begin{equation} \label{5.14}
\overline{K}(V_{\phi_1}, V_{\phi_2}) \: = \: \int_M K(\nabla \phi_1, \nabla \phi_2) \:
\left( |\nabla \phi_1|^2 \: |\nabla \phi_2|^2 \: - \: \langle \nabla \phi_1, \nabla \phi_2 \rangle^2 \right) \:
\rho \: \dvol_M \: + \: 3 \: |T_{\phi_1 \phi_2}|^2,
\end{equation}
where $K(\nabla \phi_1, \nabla \phi_2)$ denotes the sectional curvature of the
$2$-plane spanned by $\nabla \phi_1$ and $\nabla \phi_2$.
\end{corollary}

\begin{corollary} \label{5.15} If $M$ has nonnegative sectional curvature then 
$P^\infty(M)$ has nonnegative sectional curvature.
\end{corollary}

\begin{remark} \label{5.16} One can ask whether the condition of  
$M$ having sectional curvature bounded below by $r \in \R$
implies that $P^\infty(M)$ has sectional curvature bounded below by $r$.
This is not the case unless $r = 0$. The reason is one of normalizations.
The normalizations on $\phi_1$ and $\phi_2$ are
$\int_M |\nabla \phi_1|^2 \: \rho \: \dvol_M \: = \: \int_M |\nabla \phi_2|^2 \: \rho \: \dvol_M \: = \: 1$
and $\int_M \langle \nabla \phi_1, \nabla \phi_2 \rangle \: \rho \: \dvol_M \: = \: 0$. One
cannot conclude from this that 
$
\int_M 
\left( |\nabla \phi_1|^2 \: |\nabla \phi_2|^2 \: - \: \langle \nabla \phi_1, \nabla \phi_2 \rangle^2 \right) \:
\rho \: \dvol_M
$
is $\ge 1$ or $\le 1$.

More generally, if $M$ has nonnegative sectional curvature then
$P_2(M)$ is an Alexandrov space with nonnegative curvature
\cite[Theorem A.8]{Lott-Villani}, \cite[Proposition 2.10(iv)]{Sturm}. On the other hand,
if $M$ does not have nonnegative sectional curvature then one sees by an
explicit construction that $P_2(M)$ is not an Alexandrov space with curvature
bounded below \cite[Proposition 2.10(iv)]{Sturm}.
\end{remark}

\begin{remark} \label{5.17}
The formula (\ref{5.5}) has the structure of the O'Neill formula for the sectional curvature
of the base space of a Riemannian submersion.  In the case $M = \R^n$, Otto
argued that $P^\infty(\R^n)$ is formally the quotient space of $\Diff(\R^n)$, with an $L^2$-metric,
by the subgroup that preserves a fixed volume form \cite{Otto}. As   $\Diff(\R^n)$ is formally
flat, it followed that $P^\infty(\R^n)$ formally had nonnegative sectional curvature. 
\end{remark}

\section{Poisson structure} \label{section6}

Let $M$ be a smooth connected closed manifold.  We do not give it a Riemannian metric.
In this section we describe a natural Poisson structure on 
$P^\infty(M)$ arising from a Poisson structure on $M$. If $M$ is
a symplectic manifold then we show that the symplectic leaves in $P^\infty(M)$
are orbits of the action of the group $\Ham(M)$ of Hamiltonian diffeomorphisms
acting on $P^\infty(M)$. We recover the symplectic structure on the orbits
that was considered in \cite{Ambrosio-Gangbo,Gangbo-Pacini}.

Let $M$ be a smooth manifold and let $p \in C^\infty({\wedge}^2 TM)$ be a skew bivector
field.  Given $f_1, f_2 \in C^\infty(M)$, one defines the Poisson bracket
$\{ f_1, f_2 \} \in C^\infty(M)$ by
$\{ f_1, f_2 \} \: = \: p(df_1 \otimes df_2)$. There is a skew trivector field
$\partial p \in C^\infty(\wedge^3 TM)$ so that for $f_1, f_2, f_3 \in C^\infty(M)$,
\begin{equation} \label{6.1}
(\partial p)(df_1, df_2, df_3) \: = \: \{ \{f_1, f_2 \}, f_3 \} \: + \: 
 \{ \{f_2, f_3 \}, f_1 \} \: + \: 
 \{ \{f_3, f_1 \}, f_2 \}.
\end{equation}
One says that $p$ defines a Poisson structure on $M$ if $\partial p \: = \: 0$.
We assume hereafter that $p$ is a Poisson structure on $M$.

\begin{definition} \label{6.2}
Define a skew bivector field $P \in C^\infty(\wedge^2 TP^\infty(M))$ by saying
that its Poisson bracket is 
$\{ F_{\phi_1}, F_{\phi_2} \} \: = \: F_{\{\phi_1, \phi_2\}}$, i.e.
\begin{equation} \label{6.3}
\{ F_{\phi_1}, F_{\phi_2} \} (\mu) \: = \: \int_M \{\phi_1, \phi_2\} \: d\mu
\end{equation}
for $\mu \in P^\infty(M)$.
\end{definition}

The map $\phi \rightarrow dF_\phi \Big|_{\mu}$ passes to an isomorphism
$C^\infty(M)/\R \rightarrow T^*_{\mu} P^\infty(M)$. As
the right-hand side of (\ref{6.3}) vanishes if $\phi_1$ or $\phi_2$ is constant, equation
(\ref{6.3}) does define an element of $C^\infty(\wedge^2 TP^\infty(M))$.

\begin{proposition} \label{6.4}
$P$ is a Poisson structure on $P^\infty(M)$.
\end{proposition}
\begin{proof}
It suffices to show that $\partial P$ vanishes. This follows from the equation
\begin{align} \label{6.5}
(\partial P)(dF_{\phi_1}, dF_{\phi_2}, dF_{\phi_3}) \: & = \: 
\{ \{F_{\phi_1}, F_{\phi_2} \}, F_{\phi_3} \} \: + \: 
\{ \{F_{\phi_2}, F_{\phi_3} \}, F_{\phi_1} \} \: + \: 
\{ \{F_{\phi_3}, F_{\phi_1} \}, F_{\phi_2} \}  \\ 
 &  = \: F_{\{ \{\phi_1, \phi_2 \}, \phi_3 \} \: + \: 
 \{ \{\phi_2, \phi_3 \}, \phi_1 \} \: + \: 
 \{ \{\phi_3, \phi_1 \}, \phi_2 \}} \: = \: 0. \notag
\end{align}
\end{proof}

A finite-dimensional Poisson manifold has a (possibly singular) foliation with
symplectic leaves \cite{Kirillov}. The leafwise tangent vector fields are spanned by
the vector fields $W_f$ defined by $W_f h \: = \: \{f,h\}$. The symplectic form $\Omega$
on a leaf is given by saying that 
$\Omega(W_{f}, W_{g} )\: = \: \{ f, g \}$.

Suppose now that $(M, \omega)$ is a closed $2n$-dimensional symplectic manifold. 
Let $\Ham(M)$ be the group of Hamiltonian symplectomorphisms of $M$
\cite[Chapter 3.1]{McDuff-Salamon}.
\begin{proposition} \label{6.6}
The symplectic leaves of $P^\infty(M)$ are the orbits of the action of
$\Ham(M)$ on $P^\infty(M)$. Given $\mu \in P^\infty(M)$ and 
$\phi_1, \phi_2 \in C^\infty(M)$, let $\widehat{H}_{\phi_1}, \widehat{H}_{\phi_2} \in 
T_{\mu} P^\infty(M)$ be the
infinitesimal motions of $\mu$ under the flows generated
by the Hamiltonian vector fields $H_{\phi_1}, H_{\phi_2}$ on $M$. Then
$\Omega(\widehat{H}_{\phi_1}, \widehat{H}_{\phi_2}) \: = \:
\int_M \{ \phi_1, \phi_2 \} \: d\mu$.
\end{proposition}
\begin{proof}
Write $\mu \: = \: \rho \: \omega^n$. We claim that $(W_{F_\phi} \widehat{F})(\mu) \: = \:
\frac{d}{d\epsilon} \Big|_{\epsilon = 0} \widehat{F}(\mu \: - \: \epsilon \: \{\phi,\rho\} \: \omega^n)$
for $\widehat{F} \in C^\infty(P^\infty(M))$.
To show this, it is enough to check it for each $\widehat{F} \: = \: F_{\phi^\prime}$, with
$\phi^\prime \in C^\infty(M)$. But
\begin{equation} \label{6.7}
(W_{F_\phi} {F_{\phi^\prime}})(\mu) \: = \: F_{\{\phi, \phi^\prime \}}(\mu) \: = \:
\int_M \{ \phi, \phi^\prime \} \: \rho \: \omega^n \: = \: - \:
\int_M \phi^\prime \: \{ \phi, \: \rho \} \: \omega^n,
\end{equation}
from which the claim follows. This shows that $W_{F_\phi} \: = \: \widehat{H}_\phi$.

Next, at $\mu \in P^\infty(M)$ we have
\begin{equation} \label{6.8}
\Omega(\widehat{H}_{\phi_1}, \widehat{H}_{\phi_2}) \: = \:
\Omega(W_{F_{\phi_1}}, W_{F_{\phi_2}}) \: = \:
\{ F_{\phi_1}, F_{\phi_2} \} (\mu) \: = \: 
\int_M \{ \phi_1, \phi_2 \} \: d\mu.
\end{equation}
This proves the proposition.
\end{proof}

\begin{remark} \label{6.9} As a check on Proposition \ref{6.6}, suppose that
$\phi_2 \in C^\infty(M)$ is such that $\widehat{H}_{\phi_2}$ vanishes at
$\mu \: = \: \rho \: \omega^n$. Then $\{ \phi_2, \rho \} = 0$, so by our
formula we have
\begin{equation} \label{6.10}
\Omega(\widehat{H}_{\phi_1}, \widehat{H}_{\phi_2}) \: = \:
\int_M \{ \phi_1, \phi_2 \} \: d\mu \: = \: 
\int_M \{ \phi_1, \phi_2 \} \: \rho \: \omega^n \: = \: 
\int_M  \phi_1 \: \{ \phi_2,  \rho \} \: \omega^n \: = \: 0.
\end{equation}
\end{remark}

\begin{remark} \label{6.11}
The Poisson structure on $P^\infty(M)$ is the restriction of the Poisson structure
on $(C^\infty(M))^*$ considered in 
\cite{Marsden-Weinstein,Marsden-Ratiu-Schmid-Spencer-Weinstein,Weinstein}.
Here the Poisson structure on $(C^\infty(M))^*$
comes from the general construction of a Poisson structure on the dual of a Lie
algebra, considering $C^\infty(M)$ to be a Lie algebra with respect to
the Poisson bracket on $C^\infty(M)$.
The cited papers use the Poisson structure
on $(C^\infty(M))^*$ to show that certain PDE's are Hamiltonian flows.

\end{remark}

\end{document}